\documentclass[11pt]{amsart}

\usepackage{graphicx}
\usepackage{latexsym}

\usepackage{amssymb,amsmath}
\usepackage{epstopdf}
\usepackage{hyperref}
\usepackage{color}
\usepackage{tikz}

\input glyphtounicode
\usetikzlibrary{positioning}
\usetikzlibrary{matrix}
\newtheorem{theorem}{Theorem}[section]
\newtheorem{def.}[theorem]{Definition}
\newtheorem{prop.}[theorem]{Proposition}
\newtheorem{lem.}[theorem]{Lemma}
\newtheorem{cor.}[theorem]{Corollary}
\newtheorem{conj.}[theorem]{Conjecture}
\newtheorem{Bsp.}{Example}[section]
\newtheorem{rem.}[theorem]{Remark}

\def\Hil{\mathcal{H}}
\def\K{\mathcal{K}}
\def\B{\mathcal{B}}

\def\D{\mathcal{D}}
\newcommand{\ran}{{\sf Ran}}
\newcommand{\dom}{{\sf Dom}}
\newcommand{\spa}{{\sf span}}

\newcommand{\gl}{{\mathfrak L}}
\newcommand{\LDD}{\gl(\D,\D^\times)}
\newcommand{\ad}{^{\mbox{\scriptsize $\dag$}}}

\newcommand{\identity}[1]{\mathsf{id}_{ #1 }}

\newcommand{\ip}[2]{\left\langle {#1}\left | {#2}\right.\right\rangle}

\def\io{\iota}
\def\paz{\bigstar}
\definecolor{Magenta}{rgb}{1,0,1}
\definecolor{xxl}{rgb}{0.2,0.08,0.4}

\begin{document}
\title{Frame-related sequences in chains and scales of Hilbert spaces}
\author{Peter Balazs}
\author{Giorgia Bellomonte} 
\author{Hessam Hosseinnezhad}
\address{Acoustics Research Institute, Austrian Academy of Sciences, Wohllebengasse 12-14,	A-1040 Vienna, Austria; }\email{peter.balazs@oeaw.ac.at}
\address{
	Dipartimento di Matematica e Informatica, Universit\`{a} degli Studi di Palermo,
	I-90123 Palermo, Italy; }\email{giorgia.bellomonte@unipa.it}
\address{ Department of Basic Sciences, Shahid Sattari Aeronautical University of Science and Technology, Tehran, Iran; }
\email{hosseinnezhad@ssau.ac.ir}

%\date{\today}
\begin{abstract}Frames for Hilbert spaces are interesting for mathematicians but also important for applications e.g. in signal analysis and in physics. 
	Both in mathematics and physics it is natural to consider a full scale of spaces, and not only a single one. 
	In this paper, we study how certain frame-related properties, as completeness or the property of being a (semi-)frame, of a certain sequence in one of the spaces propagate to other spaces in a scale. We link that to the properties of the respective frame-related operators, like analysis or synthesis. We start with a detailed survey of the theory of Hilbert chains. Using a canonical isomorphism the properties of frame sequences are naturally preserved between different spaces. We also show that some results can be transferred if the original sequence is considered, in particular that the upper semi-frame property is kept in larger spaces, while the lower one to smaller ones. This leads to a negative result: a sequence can never be a frame for two Hilbert spaces of the scale if the scale is non-trivial, i.e. spaces are not equal. \end{abstract}
\maketitle

\section{Introduction}

Frames have been used as a powerful alternative to Hilbert space bases, and they allow a deep theory (for an overview see \cite{caku13,ole1n,he98-1}). They  are also very important for applications, e.g. in physics \cite{antoin2,colgaz10}, signal processing \cite{nsdgt10,benli98,boelc1}, numerical treatment of operator equations \cite{dafora05,Stevenson03} and acoustics \cite{xxllabmask1,majxxl10}. 

There have been numerous generalizations of the concept of frames, see e.g. \cite{antoin2, GBell, BC, cakuli08,chst03, gavruta}, among others. \\

The  basic idea in the work by Duffin and Schaeffer \cite{duffschaef1} was to have a sequence of elements in a Hilbert space, that allow redundant and stable representations.
It is often more natural to not only consider a single Hilbert space, but a whole chain or scale of spaces \cite{Antoine1998,grhe99,raul11,xxlgrospeck19}. Therefore,  aiming  at an extension of the concept of frames to such a setting is very natural. Several approaches have already been done, see Gelfand frames \cite{dafora05}, or Riesz-like bases in rigged Hilbert spaces \cite{belltrap18}. Those concepts ``only'' deal with a triplet of spaces, while here we will work on the concept of frames and related objects on a full scale.
Applications of frames for scales of Hilbert spaces can be found in the discretization of operators (as in \cite{xxlframoper1,xxlgro14,xxlhar18,harbr08,Stevenson03}).

For localized frames \cite{xxlgro14,forngroech1} a natural scale of (Banach) spaces is associated to a frame in a central Hilbert space. For this concept the frame-related properties get naturally shared on all spaces. In this manuscript we investigate how frame-related properties are transferred in a general scale of Hilbert space, motivated by the applications using e.g. Sobolev spaces \cite{eh10}.\\

Let the increasing {\em chain of Hilbert spaces} be given \cite{Antoine1998}:
$$ ... \subseteq \Hil_2 \subseteq \Hil_1 \subseteq \Hil_0 \subseteq \Hil_{-1} \subseteq \Hil_{-2} \subseteq ... , $$
with dense inclusions and  $\Hil_{-n} = \Hil_n^\times$, where $\Hil_n^\times$  is the antidual/conjugate dual of $\Hil_n$ with respect to the inner product of $\Hil_0$ \cite{ait_book}  i.e. the space of continuous conjugate
linear functionals on the Hilbert space $\Hil_n$ \cite%[p. 254]
{Shmudg}.

If such a chain is generated by the domain of an operator \cite{antoine2009partial}, we call this a {\em scale of Hilbert spaces}. Note that, a chain of three Hilbert spaces  generates always a scale of Hilbert spaces, see Section \ref{subsec: putting all together}. For a similar definition, see the one of nested Hilbert space \cite{CIS-293226}, where only a partial order is assumed. See also the similar concept of Gelfand chain \cite{mcpi11}. 
In this paper, we will look at frame-like sequences for such a scale, and how properties of those sequence are transferred between spaces.

This paper is organized as follows: In Section \ref{sec:not0} we collect results from the literature and fix the notation. In Section \ref{sec:Hilchain0} we give a survey for chains of Hilbert spaces, following literature to some extent but describing it from a new point of view, focusing on the so-called Berezanskii isomorphisms between two arbitrary spaces of the scale. In Section \ref{propagation} we study how certain frame-related properties of a certain sequence in a space of a scale of Hilbert spaces and of operators directly linked to it propagate in the whole scale.

\section{Known facts, definitions and notation} \label{sec:not0}
Before going forth, let us introduce some notations and recall the main definitions in literature we are going to use.\\ 

Let $\Hil, \K$ be two separable infinite dimensional Hilbert spaces,  with inner 
products $\ip{\cdot}{\cdot}_\Hil$ resp. $\ip{\cdot}{\cdot}_\K$, chosen
to be linear in the first entry, and the induced norms $\|\cdot\|_\Hil$ resp. $\|\cdot\|_\K$. 
A bounded operator $A
\in\B(\Hil,\mathcal{K})$ is said to be {\em unitary} if it is isometric, i.e. $
\|Af\|_\K=\|f\|_\Hil$ for every $f\in\Hil$ and is onto,  i.e. its range $\ran(A)=
\mathcal{K}$. { If we refer to a full sequence we will denote it by a letter without an index, i.e. $c = \left(c_k\right)$.}

\subsection{Frames in Hilbert Spaces}
Let now $\Hil$ be a separable Hilbert space with inner product $\ip{\cdot}{\cdot}$ and norm $\|\cdot\|$.
A  sequence $(\psi_k)\subset \Hil$ is said to be 
\begin{itemize}
	\item {\em complete} (or {\em total}) if $\spa(\psi_k)$, the linear span of $(\psi_k)$, is dense
	in $\Hil$;
	\item  a {\em frame for $\Hil$}
	if there exist $A> 0$ and $B<\infty$ such that 
	\begin{equation}\label{eq: frame} A\|f\|^2\leq \sum \limits_{k \in \mathbb{N}} |\ip{f}{\psi_k}|^2 \leq B\|f\|^2, \qquad \forall f\in \Hil;\end{equation}
	\item a {\em Bessel sequence in $\Hil$} if 
	there exists $B>0$ such that  the upper inequality in \eqref{eq: frame} holds true.
	It is called an {\em upper semi-frame} \cite{jpaxxl09}, if the Bessel sequence is also complete (this is equivalent to $0<\sum \limits_{k \in \mathbb{N}} | \ip{f}{\psi_k} |^2 \leq B ||f||^2,\,\, \forall f \not=0$);
	\item a {\em lower semi-frame} for $\Hil$ if it satisfies the lower frame inequality in \eqref{eq: frame};
	\item a {\em Riesz basis for $\Hil$} if there exist  an orthonormal basis $(e_k)$ for $\Hil$
	and a bounded
	bijective operator $T :\Hil \to\Hil$ such that $\psi_k=Te_k$ for all $k \in \mathbb{N}$. 
\end{itemize}

There are several operators canonically associated to a sequence $\psi=(\psi_{k})$ of elements of a Hilbert space $\Hil$. The {\it analysis operator} $C_\psi :\dom(C_\psi )\subseteq \Hil\to \ell^2$ of $(\psi_{k})$ is defined by 
$$C_\psi f=(\ip{f}{\psi_{k}}),\qquad \forall f\in \dom(C_\psi), \text{ where }$$
$$\dom(C_\psi )=\left \{f\in \Hil: \sum_{k \in \mathbb{N}} |\ip{f}{\psi_{k}}|^2<\infty\right  \}.$$ 
The {\it synthesis operator} $D_\psi:\dom (D_\psi)\subseteq \ell^2 \to \Hil$ of $(\psi_{k})$ is defined on the dense domain
$$\dom(D_\psi):=\left \{ { c } \in \ell^2:\sum_{k \in \mathbb{N}} c_{k} \psi_{k} \text{ is convergent in } \Hil\right \}$$  by
$$D_\psi(c_{k})=\sum_{k \in \mathbb{N}} c_{k} \psi_{k},\qquad\forall(c_{k})\in \dom(D_\psi).$$

It is known \cite{xxlstoeant11} that we have for {\em any} sequence $C_\psi = D_\psi^*$  where, as usual, $D_\psi^*$ indicates the adjoint of the operator $D_\psi$.
Other two operators associated to a sequence $\psi=(\psi_{k})\subset\Hil$ are the {\em frame operator}\footnote{Note that for the definition of the frame operator as a potentially unbounded operator the sequence does not have to be a frame.}
$S_\psi : \dom(S_\psi)\subseteq\Hil \to \Hil$ of $\psi$

$$S_\psi f:=
\sum_{k \in \mathbb{N}}
\ip{f}{\psi_k} \psi_k$$ where
$$\dom(S_\psi) = \{f \in \Hil :
\sum_{k \in \mathbb{N}}
\ip{f}{ \psi_k} \psi_k \mbox{ converges in }\Hil\} ;$$  and the {\em Gram operator}
$ G_\psi : \dom(G_\psi)\subseteq\ell^2 \to \ell^2,$ with $$ { \left( G_\psi c \right)_k := 
	\sum_{l \in \mathbb{N}}(G_\psi)_{k,l} \cdot c_l },$$ where
{	\begin{equation*}
		\dom(G_\psi) =\{ {c} \in \ell^2 :  \sum_{l \in \mathbb{N}} (G_\psi)_{k,l}c_l\mbox{ converges }\forall k \in \mathbb{N}\mbox { and is in }
		\ell^2\}
	\end{equation*}
}

and the Gram matrix $((G_\psi)_{k,l})_{k,l}$ is defined by $(G_\psi)_{k,l} =\ip{\psi_l}{\psi_k}$, ${k,l} \in \mathbb{N}$.

If we combine those operators we end up at the following definition: 
let $\psi = (\psi_k)$ and $\phi = (\phi_k)$ be two sequences in $\Hil$. The pair $\left( \psi, \phi\right)$ is called 
\begin{itemize}
	\item a {\em reproducing pair} if the cross-frame operator defined by 
	$$ \left< S_{\psi,\phi} f,g \right> = \sum \limits_k \left< f, \psi_k \right> \left< \phi_k, g \right>$$
	is an invertible, bounded operator \cite{spexxl14}.
	\item a {\em weakly dual} pair \cite{feizim1} if it is a reproducing pair and $S_{\psi,\phi} = \identity{}$.
\end{itemize}
Clearly, for every reproducing pair, the pair $\left( \psi, S_{\psi,\phi}^{-1} \phi \right)$ is weakly dual.

\subsection{Rigged Hilbert spaces}

In a formulation using a topology viewpoint, where not all involved spaces have to be normed, see e.g. \cite{belltrap18,TraTri19}, a {\em rigged Hilbert space} (RHS) consists of a triplet $(\D
,\Hil, \D^\times)$ where $\D$ is a dense subspace of $\Hil$ endowed with
a locally convex topology $t$, finer than that induced by the Hilbert norm of $\Hil$,
and $\D^\times$ is the conjugate dual of $\D[{t}]$, endowed with the
strong topology ${t}^\times:= \beta(\D^\times, \D)$. We have \begin{equation*}\label{incl}
	\D[t] \hookrightarrow  \Hil \hookrightarrow\D^\times[t^\times],
\end{equation*} where $\hookrightarrow $ denotes a continuous
embedding and since $\D^\times$ contains a subspace that can be identified with $\Hil$, we will read \eqref{incl} as a chain of topological inclusions: $\D[t]
\subset  \Hil \subset\D^\times[t^\times]$. These identifications
imply that the sesquilinear form $B( \cdot , \cdot )$ that puts $\D$
and $\D^\times$ in duality is an extension of the inner product of
$\D$: $B(f, g) = \ip{f}{g}$, for every $f,g \in \D$ (as usual, to simplify notations we adopt the symbol $\ip{\cdot}{\cdot}$ for both of
them).

Let now $\D[t] \subset \Hil \subset \D^\times[t^\times]$ be a rigged
Hilbert space,
and let $\LDD$ denote the vector space of all continuous linear maps from $\D[t]$ into  $\D^\times[t^\times]$. {If $\D[t]$ is barreled (e.g.~reflexive)}, an involution $X \mapsto X\ad$ can be introduced in $\LDD$  by the equality
\begin{equation}
	\label{eq: X dag}\ip{X\ad \eta}{ \xi} = \overline{\ip{X\xi}{\eta}}, \quad \forall \xi, \eta \in \D
\end{equation}  and  $\LDD$ results to be a $^\dagger$-invariant vector space. \\
We will show in Section \ref{3.3.1} how this can be expressed and extended in a Hilbert space setting.

\subsection{Scales of Hilbert Spaces} \label{sec:scalHil0}
Let $\Hil_0$ be a Hilbert space with inner product $\left< . , .\right>_0$.
Let $A$ be a self-adjoint, strictly positive, unbounded operator with domain $\dom(A)\subset\Hil_0$ and range  $\ran(A) \subset \Hil_0$, without loss of generality let $A = A^* \geq 1$. 
By assumption, this operator has dense domain, dense range and is closed. 
It has closed range, and is therefore onto. 
Its inverse $A^{-1}: \Hil_0 \rightarrow \dom(A)$ is bounded and self-adjoint.

The same argument holds for any $A^k$ and $A^{-k}$, $k\in\mathbb{N}$. 
Define, for each $k \geq 0$, the Hilbert space $\Hil_{+k}=\left\{ \dom(A^{k/2}) , \|\cdot\|_{+k} \right\}$ i.e. the domain of $A^{k/2}$ equipped with the norm $\|\cdot\|_{+k}$ induced by the inner  product
$\ip{x}{y}_{+k} := \ip{A^{k/2} x}{A^{k/2} y}_0, x,y\in \dom(A^{k/2})$.
Therefore $A^k:\Hil_{+2k} \rightarrow \Hil_0$ is bounded, even unitary. 

Define $\Hil_{-k}$ as the completion of $\Hil_0$ with respect to the
norm $\|\cdot\|_{-k}$ induced  by the inner product
$\ip{f}{g}_{-k} :=\ip{A^{-{k/2} }f}{A^{-{k/2}}g}_0, f,g\in\Hil_0$.
Clearly, $\Hil_0 \subset (\Hil_{+k})'$ densely (where we apply the Riesz isomorphism for $\Hil_0 \cong (\Hil_0)'$).
We have for $f \in \Hil_0$
$$\| f \|_{\left(\Hil_{+k}\right)'} = \sup \limits_{\|g\|_{+k}=1} | \left< f,g \right>_0 | = \sup \limits_{\|A^{k/2} g\|_{0}=1} \left| \left< f,g \right>_0 \right| =  $$
$$ = \sup \limits_{\|h\|_{0}=1} | \left< f,A^{-k/2}  h \right>_0|  = \sup \limits_{\|h\|_{0}=1} | \left< A^{-k/2}  f, h \right>_0| = \| A^{-k/2} f \|_{0}  .$$
and so the norms are equivalent. 
Therefore, $(\Hil_{+k})'=\Hil_{-k}$.
For each fixed $k > 0$ the triplet
$$\Hil_{+k} \subset \Hil_0 \subset\Hil_{-k}$$
is a rigged Hilbert space with dense inclusions.
The chain of Hilbert spaces infinite in both sides
$$...\subset \Hil_{+k}\subset ...\subset \Hil_0\subset ...\subset \Hil_{-k}\subset ... $$
is called $A$-{\em scale of Hilbert spaces} \cite{Koshmanenko-2016}. 

An important property of an $A$-scale
is the invariance of the structure of rigged triple $\mathcal{H}_{+k}\subseteq\Hil_{0}\subseteq\mathcal{H}_{-k}$
under shift along the $A$-scale, i.e. shift of the index $k$, this property has been called the \emph{first invariance principle} of the $A$-scale \cite{Koshmanenko-2016}:
for any fixed $k>0$ and an arbitrary $n$, the triple of spaces $\mathcal{H}_{n+k}\subseteq\mathcal{H}_{n}\subseteq\mathcal{H}_{n-k}$
is a new rigged Hilbert space, where in particular  the 
Hilbert space
$\Hil_{n-k}$ 
is the dual of $\Hil_{n+k}$, which becomes apparent, if we do the above construction taking $\Hil_n$ as the ``central space''.

\section{Hilbert chains} \label{sec:Hilchain0}
In this section we include a detailed introduction to Hilbert chains, to some extent following \cite{be68}, albeit extended and reformulated. 
We have done this to make the manuscript self-contained, stress results not well known in the frame theory community and adapt the respective point-of-view. We have added details and covered new side-aspects. 

\subsection{Hilbert Triplets} \label{sec:isochain1}
This section deals with canonical Hilbert triplets and how one can define operators between them. 

Let $\mathcal{H}_0$ be a separable, infinite dimensional Hilbert space with scalar product $\ip{\cdot}{\cdot}_{0}$ and norm $\|\cdot\|_0$.
Fix any $i\in\mathbb{N}\setminus\{0\}$,  and let $\Hil_{+i}$ be a dense subspace of $\Hil_{0}$, $\Hil_{+i} \subseteq \Hil_0$, complete w.r.to the norm $\|\cdot\|_{+i}$ induced by the inner product $\ip{\cdot}{\cdot}_{+i}$
with \begin{equation}\label{eq: ineq norms}
	\|x\|_0\leq\|x\|_{+i},\quad\forall x\in \Hil_{+i}.
\end{equation} 
Let $\iota_{+i,0}$ be the inclusion of $\Hil_{+i}$ into $\Hil_0$, which by \eqref{eq: ineq norms} is bounded. 
We have
\begin{equation} \label{eq: dual inner prod} \ip{f}{x}_0 = \ip{f}{\iota_{+i,0} x}_0 = \ip{\iota_{+i,0}^* f}{x}_{+i},\quad \forall f\in\Hil_0, x\in\Hil_{+i}.
\end{equation}
Define
\begin{equation} \label{eq:inclusionadj1}
	\iota_{0,+i} := \iota_{+i,0}^* : \Hil_0 \rightarrow \Hil_{+i},
\end{equation}
where we consider the Hilbert space adjoints, i.e. apply the Riesz isomorphisms for both Hilbert spaces. 

We have $\|\iota_{+i,0}\|\leq1$ and $\iota_{+i,0}$ is injective (one-to-one) and has dense range, but cannot be onto.
Therefore $\iota_{0,+i} := \iota_{+i,0}^*$ has the same properties. 

On the other hand, let us note that, naturally,
\begin{equation*}%\label{eq:ipip}
	\ip{x}{y}_{+i} \not= \ip{\iota_{+i,0} x}{\iota_{+i,0} y}_0,\quad x,y\in\Hil_{+i}
\end{equation*}
because else \eqref{eq: ineq norms} would be an equality and the  two spaces
%whole scale (here we have not taken a scale, yet)
would collapse to only one space.

Let us define now the scalar product \begin{equation}\label{eq: dual inner prod -i} 
	\ip{f}{g}_{-i}:=\ip{\iota_{0,+i} f}{g}_0=\ip{\iota_{0,+i} f}{\iota_{0,+i}g}_{+i},\quad f,g\in\Hil_0
\end{equation} and consider the completion $\Hil_{-i}$ of $\Hil_{0}$ with respect to the norm $\|\cdot\|_{-i}$ induced by the inner product defined in \eqref{eq: dual inner prod -i}. Then, from $ \|\iota_{0,+i}\|\leq1$ it follows that  $	\|f\|_{-i}\leq\|f\|_{0}$ $\forall f\in \Hil_{0}$, and we have:
$$\Hil_{+i}\subset\Hil_0\subset\Hil_{-i},$$
with dense inclusions,  and 
\begin{equation*}%\label{eq: boundendness o f i0-1}
	\|x\|_{-i}\leq\|x\|_0\leq\|x\|_{+i},\quad\forall x\in \Hil_{+i}.
\end{equation*} 
Since $\iota_{0,+i}$ maps a dense subspace of $\Hil_{-i}$ (i.e. $\Hil_0$) into $\Hil_{+i}$ we consider the (unique!) bounded extension $I_{-i,+i} = \overline{\iota_{0,+i}} \in \B(\Hil_{-i},\Hil_{+i})$.

As $\iota_{0,+i}$ has dense range, $I_{-i,+i}$ is onto (surjective). By \eqref{eq: dual inner prod -i}, $I_{-i,+i}$ is an isometry,
$$\ip{\alpha}{ \beta}_{-i} = \ip{ I_{-i,+i}\alpha}{I_{-i,+i}\beta}_{+i}, 
\quad \alpha, \beta\in \Hil_{-i}.$$

Therefore, $I_{-i,+i}$ is a unitary operator, hence  
$$I_{+i,-i} := I_{-i,+i}^*=I_{-i,+i}^{-1}\in\B(\Hil_{+i},\Hil_{-i})$$
is a unitary operator too. 
The operators $I_{-i,+i}\in\B(\Hil_{-i},\Hil_{+i})$ and $I_{+i,-i}\in\B(\Hil_{+i},\Hil_{-i})$ are called {\em  Berezanskii canonical isomorphisms} (see e.g. \cite{Koshmanenko-2016}). This is a particular instance of the Riesz isomorphism, if a particular duality (with the pivot space $\Hil_0$) is chosen.

\subsection{Duality by Pivot Spaces}

Now we show that $\Hil_0$ can be considered a pivot space of $\Hil_{-i}$ and $\Hil_{+i}$ in the sense that the scalar product $\ip{\alpha}{x}_0$ defines a duality relation. Consider the bilinear form $b(f,x) = \ip{f}{x}_0$ defined on $\Hil_0 \times \Hil_{+i}$. Extend it by continuity to the bilinear form 
$$B_{-i,+i}:(\alpha,x)\in\Hil_{-i}\times\Hil_{+i}\to \mathbb{C}.$$

For $f\in\Hil_{0}$ and $x\in\Hil_{+i}$ we have 
\begin{eqnarray*}%\label{eq: inner prod 0}
	|\ip{f}{x}_{0}|&=&|\ip{\iota_{0,+i}f}{x}_{+i}|=|\ip{I_{-i,+i}f}{x}_{+i}|\\\nonumber &\leq&\|I_{-i,+i}f\|_{+i}\|x\|_{+i}=\|f\|_{-i}\|x\|_{+i}
\end{eqnarray*} 
By a limit argument we obtain
\begin{equation} \label{starstar}
	|B_{-i,+i} (\alpha,x)| \leq \|\alpha\|_{-i}\|x\|_{+i},
\end{equation}
for $\alpha\in\Hil_{-i}, x\in\Hil_{+i}$. Therefore $B_{-i,+i}$ is continuous.

\begin{rem.} \label{dual1}
	The form $B_{-i,+i}(\cdot,\cdot)$ that puts $\Hil_{+i}$ and $\Hil_{-i}$ in duality is an extension of the inner product $\ip{\cdot}{\cdot}_0$ and we will use the latter symbol for both of them.  
\end{rem.}

By a limit argument (and the conjugate symmetry of any inner product), we obtain for $\alpha, \beta\in\Hil_{-i}, x,y\in\Hil_{+i}$:

\begin{equation}\label{eq: inner prod 0 +i}
	\ip{ \alpha}{x}_{0} = \ip{ I_{-i,+i}\alpha}{ x}_{+i} \text{  and }  \ip{ x}{ y}_{+i} = \ip{ I_{+i,-i}x}{ y}_{0}
\end{equation}
and
\begin{equation*}\label{eq: relations inner pro}
	\ip{\alpha}{\beta}_{-i} = \ip{ I_{-i,+i}\alpha}{ \beta}_{0} = \ip{ \alpha}{ I_{-i,+i}\beta}_{0} =
	\ip{ I_{-i,+i}\alpha}{ I_{-i,+i}\beta}_{+i}.	
\end{equation*}

By Remark \ref{dual1} we see that $\alpha \in \Hil_{-i}$ is in $(\Hil_{+i})'$ with the same norm. On the other hand, any functional $L$ on $\Hil_{+i}$ can be represented by a $y \in \Hil_{+i}$, i.e.
$$L(x) = \ip{y}{x}_{+i} = \ip{I_{+i,-i}\,y}{x}_{0}.$$
As $I_{+i,-i}$ is an isometry we have shown:

\begin{rem.}
	This construction corresponds to considering the dual pair $(\Hil_{+i}, \Hil_0)$ using $\ip{\cdot}{\cdot}_0$ and choosing $\Hil_{-i} = \overline{\Hil_0}^{\sigma(\Hil_{+i},\Hil_0)}$ using the weak dual topology \cite{conw1}. $\Hil_{-i}$ is therefore a representation of the dual space of $\Hil_{+i}$.
\end{rem.}
Note that this is not a pure abstract baublery, but important for concrete spaces. While it is true that ``everything is isomorphic anyway'', this isomorphisms {\em cannot} be thought of as equality if concrete choices for $\Hil_{+i}$ and $\Hil_{-i}$ are worked with, for example in the setting of Sobolev spaces \cite{harbr08}.
See \cite{xxlgro14} for an application of this topology to the more general setting of coorbit spaces of localized frames. 

\begin{rem.}\label{rem: link inclusion Berez.is} 
	One might expect that there is a clear link of $I_{+i,-i}$ to the inclusion $\iota_{+i,-i}$. But clearly, the relation cannot be trivial, e.g. an extension, as the former operator is bijective, the latter is injective and has dense range and both are bounded on all of $\Hil_{+i}$.

	At the same time, it could be interesting to look at the role of the inclusion of $\Hil_0$ into $\Hil_{-i}$, i.e. $\iota_{0,-i}$ in this setting. 
	But again, it cannot be trivially linked to $I_{-i,+i}$, as e.g. $\iota_{0,-i}$ cannot be $\overline{{I_{-i,+i}}_{|_{\Hil_0}}}$ because else it would be onto. For more on that, see Remark \ref{rem:moreincluse}.
\end{rem.}

\subsection{Hilbert Chains}

Now let us add a step more. Consider a dense subspace $\Hil_{+j}$ of the Hilbert space $\Hil_{+i}$,
with $j> i$ ($i,j\in\mathbb{N}$), complete w.r. to the norm $\|\cdot\|_{+j}$,  induced on $\Hil_{+j}$ by the inner product $\ip{\cdot}{\cdot}_{+j}$,
such that $\|x\|_{+i}\leq\|x\|_{+j}$ with $x\in\Hil_{+j}$. Then  $\Hil_{+j}$ is a dense subspace of $\Hil_{0}$, too, and $\|x\|_{0}\leq\|x\|_{+j}$ with $x\in\Hil_{+j}$. If we consider the completion $\Hil_{-j}$ of $\Hil_{0}$ with respect to the norm $\|\cdot\|_{-j}$ defined like in \eqref{eq: dual inner prod -i} (with $i=j$ and  $\iota_{0,+j}:\Hil_{0}\to\Hil_{+j}\subset\Hil_{+i}$, defined like in \eqref{eq:inclusionadj1} as the adjoint of the inclusion), we have for every $ x\in \Hil_{+j}$, and every $j>i$, $i\in\mathbb{N}$:
$$\Hil_{+j}\subset\Hil_{+i}\subset\Hil_0\subset\Hil_{-i}\subset\Hil_{-j},$$   where every inclusion is dense.
Indeed, let us prove that $ \Hil_{-i}\subset\Hil_{-j}$ with $j>i$. For $f\in \Hil_{0}$ we have from \eqref{eq: dual inner prod -i} and \eqref{starstar}
$$\|f\|_{-j}^2=\ip{\iota_{0,+j} f}{f}_0\leq\|\iota_{0,+j}f\|_{+i}\|f\|_{-i} =  \|I_{-j, +j}f\|_{+i}\|f\|_{-i}  \leq $$
$$ \leq \|I_{-j, +j}f\|_{+j}\|f\|_{-i}=\|f\|_{-j}\|f\|_{-i},$$ 
for every $f\in\Hil_{0}$. 

Because  $ \Hil_{-i}\subset\Hil_{-j}$ are the closures of $\Hil_{0}$ with respect to those norms we have the dense inclusions and $$\|x\|_{-j}\leq\|x\|_{-i}\leq\|x\|_0\leq\|x\|_{+i}\leq\|x\|_{+j},\quad\forall x\in\Hil_{+j}.$$

Since we have $\io_{+j,0} = \io_{+i, 0} \io_{+j, +i}$ we also have $\io_{0,+j} = \io_{+i,+j} \io_{0,+i}$.
Again note that (in the non-trivial case) we have that $\iota_{-i, -j} I_{+i, -i} \iota_{+j, +i} \not = I_{+j, -j}$.
But, naturally, we have $\iota_{k,i} = {\iota_{j,i}}_{_{|_{\Hil_k}}}$ for $k \ge j \ge i$.

Note that we can change the role of the central pivot space:
for every $r<p$, let \begin{equation}\label{eq: iota rp}\iota_{r,p}:\Hil_r\to\Hil_p\end{equation} be the adjoint of the inclusion of the Hilbert space
$\Hil_p$ into the space $\Hil_{r}$. As it is known \begin{equation}\label{eq: ineq norms2}
	\|x\|_r\leq\|x\|_p,\quad\forall x\in \Hil_p.
\end{equation}   Consider the bilinear form $$B_{r,p}:(f,x)\in\Hil_{r}\times\Hil_p\to\ip{f}{x}_r\in\mathbb{C},$$
then by \eqref{eq: ineq norms2} it is continuous for $f\in\Hil_{r}$ and for $x\in\Hil_p$, hence it can be represented as a scalar product  both in $\Hil_r$ and in $\Hil_p$:
\begin{equation*}%\label{eq: inner products2}
	\ip{f}{x}_r=B_{r,p}(f,x)=\ip{\iota_{r,p}f}{x}_p, \quad f\in\Hil_{r}, x\in\Hil_p
\end{equation*} where $\iota_{r,p}:\Hil_r\to\Hil_p$ is a bounded operator in $\B(\Hil_r,\Hil_p)$.

\begin{rem.} \label{rem:34}
	Let $\{\Hil_{r}; {r} \in \mathbb{Z}\}$ be the family of the Hilbert spaces of a Hilbert chain.
	For every $r, p \in \mathbb{Z}$, with $p\geq {r}$ the maps $\io_{r,p}:\Hil_{r} \to \Hil_p$ are injective, such that $\|\io_{r,p} x\|_p \leq \|x\|_{r}$, $\forall x \in \Hil_{r}$, $\io_{r,r}$ is  the identity of $\Hil_{r}$ and $\io_{r,n}=\io_{p,n}\io_{r,p}$, ${r}\leq p\leq n$. Hence, the family $\{\Hil_r, \io_{r,p}, r, p \in \mathbb{Z}, p\geq r\}$ is a directed contractive system of Hilbert spaces, see \cite{BT2011}. It produces two spaces $\D$ and $\D^\times$ in conjugate duality: the latter is obtained as the inductive limit of the system, whereas the former is proved to be the projective limit of the spaces $\Hil_p$'s, with respect to the adjoint maps of the $\io_{r,p}$'s.
\end{rem.}

We will use this in the  $A$-scale framework and extend the notion of  Berezanskii canonical isomorphism for  each ordered pair of spaces $\Hil_{k}$ and $\Hil_{l}$, %$k>l$ and 
both $k,l\neq0$  by defining unitary operators as $I_{k,l}:\Hil_{k}\rightarrow\Hil_{l}$
%and its inverse $I^{-1}_{k,l}:\Hil_{+l}\rightarrow\Hil_{+k}$ 
by means of fractional powers of the operator $A$, its closure and their inverses.

\subsubsection{Different Adjoints}%{\color{Magenta}May be this subsection could be shifted above}
\label{3.3.1}
Before going forth, we give some thoughts about the adjoints, as we consider different spaces here, and ``the adjoint'' depends on which spaces are considered. \\

Consider a chain of Hilbert spaces as before. Fix any $j,i\in\mathbb{N}\setminus\{0\}$. 
For $A : \Hil_{+i}\to\Hil_{+j} $ we have already indicated by  $A^*$ the Hilbertian adjoint $A^*:\Hil_{+j}\to\Hil_{+i}$, i.e. $$\ip{ A x }{ y }_{+j}= \ip{ x}{A^* y }_{+i},\quad\forall x\in \Hil_{+i},\forall y\in\Hil_{+j}.$$

This adjoint is not always useful, as two Riesz isomorphisms are applied, which is not compatible with a structure where we consider Hilbert spaces included in each other, i.e. $ \Hil_{+j} \subseteq \Hil_{+i} \subseteq \Hil_{0} \subseteq \Hil_{-i} \subseteq  \Hil_{-j}$ but distinguished from each other. (Note that even though all Hilbert spaces are isomorphic to each other, still it makes a lot of sense to treat them differently, e.g. $L^2(\mathbb{R})$ and $\ell^2(\mathbb{N})$). 
By identifying the dual of $\Hil_{-j}$ with $\Hil_{+j}$ the whole scale would collapse.  Therefore, taking $i,j\geq0$, we now consider another adjoint $A^{\paz}\in\B(\Hil_{-j},\Hil_{-i})$ for $ A\in\B(\Hil_{+i},\Hil_{+j})$, the ``pivot adjoint'' defined as follows
\begin{equation}\label{eq: relation pivot adj n,m >0}
	\ip{  \alpha}{ A x }_{0} = \ip{ A^{\paz} \alpha}{ x }_{0} \quad\forall  \alpha\in \Hil_{-j},x\in\Hil_{+i}.
\end{equation}
Here  the number in the circle indicates the subscript of the pivot Hilbert space $\Hil_0$ w.r. to which the dual is taken. 
We can, very naturally, define the adjoint for any pivot space $\Hil_{p} \not = \Hil_{0}$.  

The same construction is possible for $A\in\B(\Hil_{-j},\Hil_{-i})$, $\B(\Hil_{-j},\Hil_{+i})$ or $\B(\Hil_{+i},\Hil_{-j})$.

We highlight that the notion of pivot adjoint is a generalization of the involuted $A^\dagger$ recalled in \eqref{eq: X dag}.

Moreover, similar to  \eqref{eq: dual inner prod}, if $ A\in\B(\Hil_{+i},\Hil_{+j})$, we have for $\alpha\in \Hil_{-j}, x\in\Hil_{+i}$
\begin{eqnarray*}%\label{eq: rel A circle 0 and A^*}
	\ip{  \alpha}{ A x }_{0} &=& \ip{I_{-j,+j}\alpha}{Ax}_{+j}=\ip{ A^*I_{-j,+j}\alpha}{x}_{+i}\\
	&=&\nonumber =\ip{I_{+i,-i}A^*I_{-j,+j}\alpha}{x}_0=\\&=&\nonumber\ip{ A^{\paz} \alpha}{ x }_{0}
\end{eqnarray*}
therefore we deduce that 
\begin{equation} \label{ggg1}
	A^{\paz}=(I_{-i,+i})^{-1}A^*I_{-j,+j}=I_{+i,-i}A^*I_{-j,+j} 
\end{equation}
and 
$$ \|A^{\paz}\|=\|A^*\|=\|A\|.$$
\begin{rem.}
	In general, if $p,q\in\mathbb{Z}$ and $A\in\B(\Hil_p,\Hil_q)$, then $A^*\in\B(\Hil_q,\Hil_p)$ and \begin{equation}\label{ggg2}
		A^{\paz}=I_{p,-p}A^*I_{-q,q}\in\B(\Hil_{-q},\Hil_{-p}).
	\end{equation}
\end{rem.}

\begin{rem.} Because $A^{\paz}$ is an adjoint, two properties follow immediately:
	\begin{enumerate}
		\item 	The double pivot adjoint of an operator $A\in\B(\Hil_p,\Hil_q)$,  $p,q\in\mathbb{Z}$, is  $$A^{\paz\paz}=(A^{\paz})^{\paz}=A\quad\mbox{ and }\quad\|A^{\paz\paz}\|=\|A\|.$$
		\item Let $A\in\B(\Hil_0,\Hil_p)$ and $B\in\B(\Hil_m,\Hil_0)$, $p,m\in\mathbb{N}$, then
		$(AB)^{\paz} = B^{\paz} A^{\paz}.$ Indeed, if $x\in\Hil_{m}$ and $y\in\Hil_{p}$, by \eqref{eq: relation pivot adj n,m >0} $$\ip{AB x}{y}_0= \ip{ Bx }{ A^{\paz} y}_{0}=\ip{ x }{B^{\paz} A^{\paz} y}_{0}.$$ 
	\end{enumerate} 
\end{rem.}

%%%%%%%%%%%%%%%%

As a side remark to \eqref{ggg1} note that $I_{+i,-i}$ is a unitary operator with respect to the pair $(\Hil_{+i},\Hil_{-i})$, but is selfadjoint with respect to $\Hil_0$.

Note that 
$A^{\paz}$ also corresponds to the Banach space adjoint of $A : \Hil_{+i} \rightarrow \Hil_{+j}$ where the Riesz isomorphism at the pivot space level is considered to be an equality, i.e. $\Hil_0 = \Hil_0^\times$.
It maps $\Hil_{-j}$ into $\Hil_{-i}$.

\begin{rem.}\label{rem:moreincluse} 
	Define
	
	\begin{equation*}\label{eq: iota 0,-i}
		\iota_{0,-i} := \iota_{+i,0}^{\paz}: \Hil_0 \rightarrow \Hil_{-i} ,
	\end{equation*}
	we have $\|\iota_{0,-i}\|\leq1$. As we have  seen in \eqref{ggg2}
	\begin{equation*}
		\iota_{0,-i}=(I_{-i,+i})^{-1}\io_{+i,0}^*I_{0,0}=I_{+i,-i}\io_{0,+i}.      
	\end{equation*}The operator $\io_{0,-i}$ is effectively the inclusion of $\Hil_0$ in $\Hil_{-i}$ indeed, by \eqref{eq: inner prod 0 +i},
	$$\ip{f}{x}_{0} = \ip{\io_{0,+i} f}{x}_{+i}=\ip{I_{+i,-i}\iota_{0,+i}f}{x}_{0},\quad f\in\Hil_0,x\in\Hil_{+i}.$$ 
\end{rem.}

\subsubsection{Putting all together}\label{subsec: putting all together}

Let us  go back to the operator $\iota_{0,+1}$; it acts continuously from $\Hil_{0}$ to $\Hil_{+1}$.
Since $\Hil_{+1}\subseteq\Hil_{0}$, this operator may be considered as acting in $\Hil_{0}$; denote it as the operator $\hat{I}: \Hil_0 \rightarrow \Hil_0$. We have 
\begin{equation}\label{eq:hat{I}}
	\hat{I} = \iota_{+1,0} \iota_{0,+1}.
\end{equation}

The operator $\hat{I}$ is continuous with bound less or equal to $1$,				positive and invertible onto $\ran(\hat{I})$, with $\ran(\hat{I}) = \dom(\hat{I}^{-1})$  dense in $\mathcal{H}_0$.

It is easy to see that $\hat{I}^{-1}$ is also self-adjoint and positive in $\mathcal{H}_0$,
(later it will be clear that $ \ran(\hat{I}) = \dom(\hat{I}^{-1})=\Hil_{+2}$, and $\hat{I}^{-1}=I_{+2,0}:\Hil_{+2}\to\Hil_0$, see \eqref{eq: inner product H_n}).

The operator $\hat{I}^{-1}$ is densely defined as an operator in $\Hil_0$, has closed range and is one-to-one. Because it has a bounded inverse, it is closed. It is positive and self-adjoint and we have 
${\left( \hat{I}^{-1} \right) }^{1/2} = {\left({\hat{I}}^{1/2}\right)}^{-1} =: \hat{I}^{-1/2}$.
The properties of this operator can be summarized by
\begin{theorem}\cite[Theorem I.1.1]{be68}\label{thm: berenz}
	Consider the operator $F = (\hat{I}^{-1})^{1/2}$ in the space $\Hil_{0}$.
	It is a positive self-adjoint operator for which $\dom(F) = \Hil_{+1}$ and
	$\ran(F) = \Hil_{0}$. This operator acts isometrically from $\Hil_{+1}$ onto $\Hil_{0}$:
	%{@Peter}
	%$(i^{-1}f,g)=<Î^{-1},Î^{f},Î^{g}> = <f*,I^{g'}=<Îf'},g'>=$
\begin{equation*}
	\ip{ x}{ y}_{+1} = \ip{ F x}{ F y}_{0}~~~(x, y\in\Hil_{+1}).
\end{equation*}
Consider $F$ as an operator acting from $\Hil_{0}$ to $\Hil_{-1}$ and
form the closure by continuity, denote this operator by $\mathbf{F}$.
$\mathbf{F}$ acts isometrically from the whole $\mathcal{H}_0$ onto $\Hil_{-1}$:
\begin{equation*}
	\ip{ f}{ g}_{0} = \ip{\mathbf{F}f}{ \mathbf{F}g}_{-1}~~~(f, g\in\Hil_{0}),
\end{equation*}
and moreover $I_{+1,-1}=I_{-1,+1}^{-1} = \mathbf{F}F$. The relation
\begin{equation*}%\label{eq: pivot F_1}
	\ip{ f}{ Fx}_{0} = \ip{\mathbf{F}f}{ x}_{0}~~~(f\in\Hil_{0}, x\in\Hil_{+1}),
\end{equation*}
holds, therefore $\mathbf{F}=F^{\paz}$.
\end{theorem}

Let us write $B = F^{-1}$ and $\mathbf{B} = \mathbf{F}^{-1}$. 

From the factorization of $I_{-1,+1}$ it follows immediately that $$I_{+1,-1}=I_{-1,+1}^{-1} = \mathbf{B}^{-1}B^{-1}=\mathbf{F}F= F^{\paz}F.$$  thus obtaining a factorization of $I_{-1,+1} $ in terms of isometric operators. 
\\

Using the Hilbert adjoint we can write
$$ B = F^*  \text{ and } \mathbf{B} = \mathbf{F}^*.$$

We have the following scheme

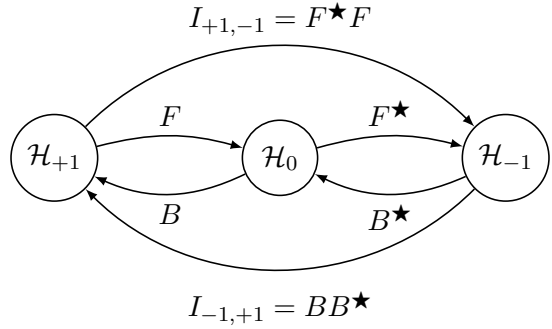
\begin{center}
\begin{tikzpicture}[-latex ,auto ,node distance=2cm and 3cm,on grid ,
	semithick ,
	state/.style={circle ,top color=white ,bottom color=white ,
		draw , black  , text=black ,minimum width=1cm}]
	\node[state] (0) {$\Hil_0$};
	\node[state] (+1) [ left=of 0] {$\Hil_{+1}$};
	\node[state] (-1) [ right=of 0] {$\Hil_{-1}$};
	\path (0) edge [bend left=25] node[below=0cm]{$B$} (+1);
	\path (+1) edge [bend right=-15 ] node[above=0cm]{$F$} (0);
	\path (+1) edge [bend left=45 ] node[above]{$I_{+1,-1}=F^{\paz}F$} (-1);
	\path (-1) edge [bend left=45 ] node[below=0.15cm]{$I_{-1,+1}=B B^{\paz}$} (+1);
	\path (0)  edge [bend left=15] node[above=0cm]{$F^{\paz}$} (-1);
	\path (-1) edge [bend right=-25 ] node[below=0cm]{$B^{\paz}$} (0);
\end{tikzpicture}
\end{center}

There is a well-known connection between infinite chains of Hilbert spaces  and positive
self-adjoint operators in $\Hil_0$, if $\Hil_0$ is the central space, see Section \ref{sec:scalHil0}.  

\begin{rem.} 
Not every  infinite chain of Hilbert spaces is a scale, see \cite[pag 161]{antoine2009partial}. 

However, see e.g. \cite[Theorem I.1.1]{be68}, given a chain of Hilbert spaces there exists a positive
self-adjoint operator $A = A^*$ on $\Hil_0$ such that the triple $\Hil_{+i}\subset\Hil_0\subset\Hil_{-i}$ of the chain coincides with the ``central'' part of the $A$-scale of Hilbert spaces generated by $A$. 
\end{rem.}
We now proceed as in Section \ref{sec:scalHil0}.\\ 

Consider the $A$-scale of Hilbert spaces where $A = F^2= \hat{I}^{-1}$.
Since $F$ is self-adjoint, then all its powers $F^n$, $n>0$ exist. For every $n>0$, the domain $\dom(F^n)$ can be made into a Hilbert space by setting
\begin{equation}\label{eq: inner product H_n}
\ip{x}{y}_{+n}=\ip{F^nx}{F^ny}_0,\qquad x,y\in\dom(F^n).
\end{equation} 
$F$ is injective, hence the following norm is induced:
$$\|x\|_{+n}=(\ip{x}{x}_{+n})^{1/2}=(\ip{F^nx}{F^nx}_0)^{1/2}=\|F^nx\|_0,\quad x\in\dom(F^n).$$
Put $\Hil_{+n}
=\dom(F^n)$, $n\geq0$. 
As usual, for every $n>0$ we can consider the Hilbert spaces $\Hil_{-n}$ the conjugate dual of $\Hil_{+n}$ with
respect to the inner product of $\Hil_0$. We get like e.g. in \cite[Example 10.1.1]{ait_book} an $A$-scale of Hilbert spaces.

We have also that for $B = F^{-1}$
\begin{equation}\label{eq: inverse}
\ip{f}{g}_{0}=\ip{B^nf}{B^ng}_{+n},\qquad f,g\in\Hil_0.
\end{equation}

By \begin{equation*}
\ip{ f}{ g}_{0} = \ip{F^{\paz} f}{ F^{\paz} g}_{-1},~~~f, g\in\Hil_{0}
\end{equation*}we deduce

\begin{equation*}
\ip{\alpha}{ \beta}_{-1}=\ip{ B^{\paz} \alpha}{ B^{\paz} \beta}_{0}, ~~~\alpha, \beta\in\Hil_{-1},
\end{equation*}

\begin{equation*}
\ip{ x}{ y}_{+1} = \ip{ Fx}{ Fy}_{0}=
\ip{F^{\paz} Fx}{ F^{\paz}  Fy}_{-1},~~~\forall x, y\in\Hil_{+1},
\end{equation*}

and 	\begin{equation*}
\ip{\alpha}{ \beta}_{-n}=\ip{ (B^{\paz}) ^n\alpha}{ (B^{\paz}) ^n\beta}_{0}, ~~~\alpha, \beta\in\Hil_{-n},\quad n>0.
\end{equation*}

\begin{rem.}\label{rem: factorization}We already know that $I_{+1,-1}=F^{\paz}F=(B^{\paz})^{-1} F$. For $0\leq r\leq p$, put $I_{p,r}=B^rF^p$. Let $p,r<0$,  put $F^p:=(F^{\paz})^p=(B^{\paz})^{-p}$ and  $B^r:=(B^{\paz})^r=(F^{\paz})^{-r} $.  
With this convention, whatever $p,r$ are in $\mathbb{Z}$, we can decompose \begin{equation*} %\label{eq:IPR1}
	I_{p,r}=B^rF^p:\Hil_p\to\Hil_r.
\end{equation*}
All these operators are unitary because they are products of unitary operators. We highlight however that they are not  Berezanskii isomorphisms unless $(p-r)/2\in\mathbb{Z}$, in fact, only in this case $\Hil_p$ and $\Hil_r$ are extreme spaces of a rigged Hilbert space.
Furthermore,   being $(B^rF^p )^*=(F^p)^*(B^r)^*=B^pF^r$, 
for $0\leq r\leq p$ by \eqref{eq: inner product H_n} and \eqref{eq: inverse} we have \begin{equation*}%\label{eq: rel betw norms} 
	\ip{f}{g}_p=\ip{F^p f}{F^p g}_0=\ip{B^rF^p f}{B^rF^p g}_r=\ip{I_{p,r}f}{I_{p,r}g}_r,\,\mbox{ with } f,g\in\Hil_p\end{equation*}
similarly, $$\ip{f}{g}_r=\ip{B^pF^r f}{B^pF^r g}_p=\ip{I_{r,p}f}{I_{r,p}g}_p,\quad\mbox{ with }f,g\in\Hil_r.$$ If now $ f\in\Hil_p, g\in\Hil_r$ \begin{equation*}\begin{split}
		\ip{I_{p,r}f}{g}_r&=\ip{B^rF^pf}{g}_r=\ip{F^p f}{F^r g}_{0}=\ip{B^pF^p f}{B^pF^r g}_p\\&=\ip{f}{B^pF^r g}_p=\ip{ f}{ I_{r,p}g}_p\end{split}\end{equation*}hence, predictably, $I_{r,p}=I_{p,r}^*$.\\
With the due changes, the same results are obtained for any $p,r\in\mathbb{Z}$.
\end{rem.} 

\subsection{Generator of a scale and shifting of the central space}

Can we shift the central space? I.e. how does the dual of $\Hil_{n_0+n}$ in respect to $\Hil_{n_0}$ look? 

\begin{cor.} \label{sec:shiftcentral1}  Let $r < p$ and let ${\Hil_p^\times}^{(r)}$ be the dual of $\Hil_p$ with respect to the topology of $\Hil_r$ (i.e. consider the triplet $\Hil_p \subseteq \Hil_r \subseteq {\Hil_p^\times}^{(r)}$).

Then ${\Hil_p^\times}^{(r)} = \Hil_{2r - p}$.
\end{cor.}
\begin{proof}
Let us consider an $A$-scale Hilbert spaces with $A = (\hat{I}^{-1})$, with $\hat{I}$ defined as in \eqref{eq:hat{I}}, so
by the first  invariant principle of $A$-scales, see the Section \ref{sec:isochain1}, for $k=p-r$ the triple
$$\mathcal{H}_{r+k} = \mathcal{H}_{p}\subseteq\mathcal{H}_r\subseteq\mathcal{H}_{2r-p} = \mathcal{H}_{r-k},$$
forms a chain of Hilbert scale.
\end{proof}

In other words: the dual of $\Hil_{n_0+n}$ in respect to $\Hil_{n_0}$ is $\Hil_{n_0-n}$.
\begin{rem.}\label{rem: dual spaces}Clearly, if $r<p$ then the space $\Hil_r$ can be considered the dual space of $\Hil_p$ with respect to some space $\Hil_{n_0}$ with $\Hil_p\subseteq\Hil_{n_0}\subseteq\Hil_r$ if and only if $p+r$ is an even number. If this is the case, $n_0=\frac{p+r}{2}$.
\end{rem.}

Now, for a fixed $A$-scale, we wonder how the operator $A$ 
changes if we fix another Hilbert space of the scale as the ``central'' space, for example $\Hil_k$, $k\in\mathbb{Z}$. Following \cite{Koshmanenko-2016} the answer to this question is given by using what has been called 
{\em the second invariance principle of the
A-scale}. %, see \cite{Koshmanenko-2016}. 
Indeed, it results that $A$ is unitarily equivalent to its image under any
``shift'' along the $A$-scale.

Just to fix some ideas let {$k\geq2$}, consider the $A$-scale $$ ... \subseteq \Hil_{2+k} \subseteq...\subseteq \Hil_{+k} \subseteq...\subseteq\Hil_{2}\subseteq...\subseteq \Hil_{2-k} \subseteq...\subseteq \Hil_{-k} \subseteq ... , $$ and the operators $$A_{+k} := I_{2+k,k}:\Hil_{2+k}\to\Hil_{+k}\quad\mbox{ and }\quad A_{-k} := I_{2-k,-k}:\Hil_{2-k}\to\Hil_{-k},$$ with $I_{p,r}$ defined as before (in  particular
$A =A_0 = I_{+2,0}$). Then it is easy to see that \cite{Koshmanenko-2016}  $A_{+k}=A_{_{|_{\Hil_{2+k}}}}$ i.e. it is the
restriction of $A$ to $\Hil_{2+k}$, hence $A_{+k}$ is self-adjoint in $\Hil_{2+k}$. The operator $A_{-k}$ is self-adjoint on $\Hil_{2-k}$ and it results $A_{-k}=\overline{A}^{\|\cdot\|_{-k}}$ the closure of $A$ in $\Hil_{-k}$. Both operators $A_{\pm k}$, $k\in\mathbb{N}$
are unitary images of the original operator $A$ on
$\Hil_0$, in fact
$$A_{+k} = I_{0,k}A I_{2+k,2}=B^kF^0B^0F^2B^2F^{2+k}=B^kF^{2+k}$$  and $$ A_{-k} = I_{0, -k}AI_{2-k,2}=B^{-k}F^0B^0F^2B^2F^{2-k}=B^{-k}F^{2-k},$$ i.e. for the sake of brevity, for $p\in\mathbb{Z}$%\setminus\{0\}$, {\color{Magenta} in fact we can include $p=0$}
$$A_{p} = B^pF^{2+p}=I_{2+p,p}.$$
We remark that the operator $A$ is essentially self-adjoint in each space $\Hil_{-k}$,
$k \geq1$.

\section{Frame-related properties on Hilbert scales}\label{propagation}

We will look at a scale of Hilbert space, and let's fix:
$$ 				\Hil_m \subseteq \Hil_p \subseteq \Hil_r, $$
i.e. $ r \le p 
\le m
$. 

We look at a sequence  $\psi = (\psi_{k}) \subseteq \Hil_m$ 
which has a property in some of the other spaces and investigate on how these properties ``spread'' to the other spaces:
if it is complete, forms a Bessel sequence, a frame, a basis, has a Riesz property, etc.. for a $\Hil_r$ (or $\Hil_p$), what can we say about this sequence in the other spaces?

If we allow the sequences to be changed the results are trivial consequences of the following straightforward generalization of \cite[Cor. 5.3.4]{ole1n}:
\begin{cor.} \label{cor:framunitop1} Given two Hilbert spaces $\Hil$ and $\K$, if $\psi = (\psi_k)$ is a frame for $\Hil$ with frame bounds $A,B$ and $U: \Hil \rightarrow \K$ is a unitary operator, then $(U \psi_k)$ is a frame for $\K$ with the same frame bounds. 
\end{cor.}

\begin{rem.} \label{rem:unitop1}
If we have a unitary operator between two Hilbert spaces, all the frame properties naturally transfer from a Hilbert space to the other.

In the sequel we will use the unitary operators $I_{r,p}$ defined like in Remark \ref{rem: factorization}, and get easy results for $\left(I_{r,p} \psi_k \right)$. See Corollary \ref{lem: frames}.
\end{rem.}
We need some preparation before that.

\subsection{Completeness}

\begin{lem.}\label{4.3.} Let $\psi = (\psi_{k}) \subseteq \Hil_m$. Then the following statements hold.
\begin{itemize}
	\item[i)] If $(\psi_{k})$ is complete in $\Hil_p$, then it is also complete in $\Hil_r$ for $r \le p \le m$.
	\item[ii)]If $(\psi_{k})$  is complete in $\Hil_r$ $r\leq m$, then $(I_{r,p}\psi_{k})$ is a complete sequence in $\Hil_p$ for any $p$.
	%$r < p \le n$.
\end{itemize}
\end{lem.}
\begin{proof} i)
Since $\Hil_p\subseteq\Hil_r$ densely, then for each $f\in\Hil_r$ and $\epsilon>0$, there is an element $x\in\Hil_p$ with $\|f-x\|_r<\epsilon/2$.
Now, let $(\psi_{k}) \subseteq \Hil_m$ be a complete sequence in $\Hil_p$. Then there exists $y\in {\spa(\psi_{k})}$ such that
$\|x-y\|_p<\epsilon/2$. So
$$\|f-y\|_r \leq \|f-x\|_r + \|x-y\|_r \leq \|f-x\|_r + \|x-y\|_p<\epsilon.$$
It follows that $\psi = (\psi_{k})$ is a complete sequence in $\Hil_r$ for $r \le p \le m$.\\
$ii)$ {Trivial by Remark \ref{rem:unitop1}.}
\end{proof}
%In what follows the hypothesis   $(p-r)/2\in \mathbb{N}$ is unnecessary.

Later, in Lemma \ref{lem: frame propag 2}, we will see that the converse of Lemma \ref{4.3.} $(i)$ is also true.

\subsection{Unbounded Frame-related operators on Hilbert chains}\label{sec: frame-rel operators}
Let us consider an arbitrary sequence $\psi = (\psi_{k}) \subseteq \Hil_{m} {\subseteq \Hil_p}$, and as in \cite{xxlstoeant11} define the analysis operator  $C_\psi^p :\dom(C_\psi^p )\subseteq \Hil_p\to \ell^2$ of $(\psi_{k})$  by $$\dom(C_\psi^p )=\left \{f\in \Hil_p: \sum_{k} |\ip{f}{\psi_{k}}_p|^2<\infty\right  \}$$   $$C_\psi^p f:=(\ip{f}{\psi_{k}}_p),\qquad \forall f\in \dom(C_\psi^p).$$

In an analogous way we can define the synthesis
operator $$D_{\psi}^p:\dom(D_{\psi}^p)\subseteq\ell^2\rightarrow\Hil_{p}$$ associated with the sequence $\psi$ by $$\dom(D_\psi^p )=\left \{c=(c_{k})\in \ell^2: \sum_{k} c_{k}\psi_{k} \mbox{ converges in }\Hil_p\right  \}$$   $$D_\psi^p c=\sum_{k} c_{k}\psi_{k},\qquad \forall c\in \dom(D_\psi^p).$$ As it is known $\dom(D_\psi^p )$ is dense in $\ell^2$, since it contains the finite sequences which form a dense subset of $\ell^2$ and $C_\psi^p=(D_{\psi}^p)^*$, hence $C_\psi^p$ is closed. \\

Clearly, we have for $r \le p$ that
\begin{equation} \label{eq:domran1} \dom{D_\psi^p} \subseteq \dom{D_\psi^r} \mbox{, \quad }\ran{D_\psi^p} \subseteq \ran{D_\psi^r} 
\end{equation}
$$\ker{D_\psi^p} \subseteq \ker{D_\psi^r} \text{ all inclusion are dense. }$$

Clearly $D_\psi^p\subset D_\psi^r$, since their domain are such that 
$  \dom{D_\psi^p} \subseteq \dom{D_\psi^r}$ and $D_\psi^p c=D_\psi^r c$ for every $c\in \dom{D_\psi^p}$. Then, if we look at $D_\psi^p$ as an operator into a subspace of $\Hil_r$ we  can say that $C_\psi^r\subset (D_\psi^p)^*_r$ where $(D_\psi^p)^*_r$ is the adjoint of $D_\psi^p$ as on operator from $\ell^2$ into $\Hil_r$.

Furthermore,
let us consider $D_\psi^{00}$ defined as operator on $c_{00}$ the space of finite sequences. 
By the above definition we have that $D_\psi^p$ is the closure of $D_\psi^{00}$ for any $p \ge m$.  So,  $D_\psi^p = \overline{D_\psi^{00}}^{\Hil_p}$. Let us now use that $r \le p$ and so $\| \cdot \|_{r} \le  \|\cdot\|_{p}$ and therefore
$D_\psi^r = \overline{D_\psi^{00}}^{\Hil_r} =  \overline{\left(\overline{D_\psi^{00}}^{\Hil_p}\right)}^{\Hil_r} = \overline{D_\psi^p}^{\Hil_r}$.
So, in summary
\begin{equation*} %\label{eq:synthov1}
D_\psi^r = \overline{D_\psi^p}^{\Hil_r},\quad \mbox{ for }
r\leq p.\end{equation*} 

Let us also introduce the combination of those operators. Consider  the ``frame operator'' $S^p_\psi : \dom(S^p_\psi)\subseteq\Hil_p \to \Hil_p$ of $\psi$

$$S^p_\psi f:=
\sum_{k \in \mathbb{N}}
\ip{f}{\psi_k}_p \psi_k$$ where
$$\dom(S^p_\psi) = \{f \in \Hil_p :
\sum_{k \in \mathbb{N}}
\ip{f}{\psi_k}_p \psi_k \mbox{ converges in }\Hil_p\} ;$$  and the Gram operator
$ G^p_\psi : \dom(G^p_\psi)\subseteq\ell^2 \to \ell^2,$ with $$ G^p_\psi(c_k)_k:=\left(
\sum_{l \in \mathbb{N}}(G^p_\psi)_{k,l}c_l\right)_k,$$ where
\begin{eqnarray*}\dom(G^p_\psi) =\{ (c_k)_k \in \ell^2 : && \sum_{l \in \mathbb{N}} (G^p_\psi)_{k,l}c_l\mbox{ converges }\forall k \in \mathbb{N}\mbox { and }\\&& \left(
\sum_{l \in \mathbb{N}}(G^p_\psi)_{k,l}c_l\right)_k
\in \ell^2\}
\end{eqnarray*}
and the Gram matrix $((G^p_\psi)_{k,l})_{k,l}$ is defined by $(G^p_\psi)_{k,l} = \ip{\psi_l}{\psi_k}_p$, ${k,l} \in \mathbb{N}$.\\

\subsection{Frame Properties of $I_{r,p}\psi$}

By using the unitary operator introduced in Remark \ref{rem: factorization} we naturally get
\begin{lem.}\label{lem: analysis on the scale} {For a given $p\in\mathbb{Z}$ let 
	$\psi = (\psi_{k}) \subseteq \Hil_p$ be an arbitrary sequence. Then for every $r \in\mathbb{Z}$, } $C^r_{I_{p,r}(\psi)} = C^p_{\psi}  I_{r,p}$ and $C^p_{\psi} = C^r_{I_{p,r}(\psi)}  I_{p,r}$.
\end{lem.}
\begin{proof}
We have
\begin{equation*}
	\begin{split}
		\dom (C^p_{\psi}  I_{r,p})
		&= \{f\in\dom (I_{r,p}): I_{r,p}f\in\dom C^p_{\psi}\}\\
		&= \{f\in\Hil_r: (\langle I_{r,p}f, \psi_{k}\rangle_p)\in\ell^2\}\\
		&= \{f\in\Hil_r: (\langle f, I_{
			p,r}\psi_{k}\rangle_r)\in\ell^2\}\\
		&= \dom(C^r_{I_{p,r}(\psi)}).
	\end{split}
\end{equation*}
And by the same argument the operators also are the same. 

The other result is obtained by multiplying  the operator $I_{p,r}$, inverse of $ I_{r,p}$, on the right in the equality $C^r_{I_{p,r}(\psi)} = C^p_{\psi}  I_{r,p}$.
\end{proof}

An analogue result is true for the other frame-related operators. We start with a result for the synthesis operator: 
%		\\{\gbc in the previous Lemma we've picked $\psi\subset\Hil_m$. Do we want to do the same choice in the following ones (4.5, 4.6, 4.7)? If yes, let's remember to specify $m\geq p,r$} {Here, we do not need a relation between the spaces, have a look, please..}

\begin{lem.}\label{lem: synthesis on the scale}  {For a given $p\in\mathbb{Z}$ let
	$\psi = (\psi_{k}) \subseteq \Hil_p$ be an arbitrary sequence. Then for every $r\in\mathbb{Z}$},  $D^r_{I_{p,r}(\psi)} = I_{p,r} D^p_{\psi} $ and $D^p_{\psi} =  I_{r,p} D^r_{I_{p,r}(\psi)} $.
\end{lem.}
\begin{proof}
If $c\in \dom(D_{I_{p,r}(\psi)}^r)\subseteq\ell^2$ then there exists $f\in \Hil_r$ such that $D_{I_{p,r}(\psi)}^r c=f$, i.e. for every $\varepsilon>0$ there exists $n_\varepsilon\in\mathbb{N}$ such that for every $n\geq n_\varepsilon$ \begin{eqnarray*}\left\|\sum_{k=1}^n c_{k}I_{p,r}\psi_{k}-f\right\|_r&=&\left\|I_{r,p}\sum_{k=1}^n c_{k}I_{p,r}\psi_{k}-I_{r,p}f\right\|_p\\&=&\left\|\sum_{k=1}^n c_{k}\psi_{k}-I_{r,p}f\right\|_p<\varepsilon.\end{eqnarray*}
Hence $D_{\psi}^p c=I_{r,p}f$, $\dom(D_{I_{p,r}(\psi)}^r)=\dom(I_{p,r} D_{\psi}^p)$ and $D_{I_{p,r}(\psi)}^r=I_{p,r} D_{\psi}^p$.
The other result is obtained by multiplying  the operator $ I_{r,p}$,  inverse of $I_{p,r}$, on the left in the equality $D_{I_{p,r}(\psi)}^r=I_{p,r} D_{\psi}^p$.\end{proof}

\begin{lem.}\label{lem: frame op on the scale}
{For a given $p\in\mathbb{Z}$  let $\psi = (\psi_{k}) \subseteq \Hil_p$  be an arbitrary sequence. Then for every $r\in\mathbb{Z}$}, $S^p_{I_{r,p}(\psi)}=I_{r,p}S^r_\psi I_{p,r}$.
\end{lem.}
\begin{proof} It is a consequence of the previous Lemma \ref{lem: analysis on the scale} and Lemma \ref{lem: synthesis on the scale} and of $(iii)$ in \cite[Proposition 3.3]{xxlstoeant11}.
\end{proof}

\begin{lem.}\label{lem: Gram op on the scale}  {For a given $p\in\mathbb{Z}$ let $\psi = (\psi_{k}) \subseteq \Hil_p$  be an arbitrary sequence. Then $C^p_\psi D^p_\psi\subseteq G^p_\psi$ and $	G^p_\psi =G^r_{I_{p,r}(\psi)}$  for every $r\in\mathbb{Z}$.}
\end{lem.}
\begin{proof} By $(iv)$ in \cite[Proposition 3.3]{xxlstoeant11} we have that $C^p_\psi D^p_\psi\subseteq G^p_\psi$, for every $p\leq n$.
Now, recall that $\ip{\psi_l}{ \psi_k}_p=\ip{I_{p,r}\psi_l}{I_{p,r} \psi_k}_r$, we have  \begin{eqnarray*}\dom(G^p_\psi)& =&\left\{ c \in \ell^2 :  \sum_{l \in \mathbb{N}} \ip{\psi_l}{ \psi_k}_p c_l\mbox{ converges }\forall k \in \mathbb{N}\mbox { and }\right. \\& & \left. 
	\sum_k \left|\sum_{l \in \mathbb{N}}\ip{\psi_l}{ \psi_k}_p c_l\right|^2
	<\infty\right\}\\
	&=&\left\{ c \in \ell^2 :  \sum_{l \in \mathbb{N}}\ip{I_{p,r}\psi_l}{I_{p,r} \psi_k}_r c_l\mbox{ converges }\forall k \in \mathbb{N}\mbox { and }\right. \\& & \left.
	\sum_k \left|\sum_{l \in \mathbb{N}}\ip{I_{p,r}\psi_l}{I_{p,r} \psi_k}_r c_l\right|^2
	<\infty\right\}=\dom(G^r_{I_{p,r}(\psi)}).
\end{eqnarray*}
\end{proof} 

Note that all the above results lead to statements saying that if $\psi$ is a frame (Bessel sequence, Riesz basis,...) for some Hilbert space $\Hil_r$, $I_{r,p} \psi$ is one for $\Hil_{p}$. They all are trivial consequences of Remark \ref{rem:unitop1}. 

\begin{cor.}\label{lem: frames}
{For a given $p\in\mathbb{Z}$ let $\psi = (\psi_{k}) \subseteq \Hil_p$ be an arbitrary sequence. Then for any $r \in \mathbb{Z}$ the following is true. 
	\begin{enumerate}
		\item If $(\psi_{k})$ is a Bessel sequence in $\Hil_p$, then $(I_{p,r}\psi_{k})$ is a Bessel sequence in $\Hil_r$.
		\item If $(\psi_{k})$ is a semi-frame in $\Hil_p$, then $(I_{p,r}\psi_{k})$ is a semi-frame in $\Hil_r$ with the same bounds.
		\item If $(\psi_{k})$ is a frame in $\Hil_p$, then $(I_{p,r}\psi_{k})$ is a frame in $\Hil_r$ with the same bounds.
		\item If   $(\psi_{k})$  and  $(\phi_{k})\subset\Hil_p$  are a reproducing pair,  then  $(I_{p,r}\phi_{k})$ and  $(I_{p,r}\psi_{k})$ are a reproducing pair in $\Hil_r$ with the same bounds.  
		\item If $(\phi_{k})\subset\Hil_p$ is a dual sequence of $(\psi_{k})$  in $\Hil_p$, then  $(I_{p,r}\phi_{k})$ is a dual  sequence of  $(I_{p,r}\psi_{k})$ in $\Hil_r$.
		\item  If $(\psi_{k})$ is an orthonormal basis of $\Hil_p$, then $(I_{p,r}\psi_{k})$ is an orthonormal basis of $\Hil_r$.
		\item If $(\psi_{k})$ is a Riesz basis of $\Hil_p$ and $T\in\B(\Hil_p)$ is the bijective operator such that $Te_{k}=\psi_{k}$ for every $k$ with $\{e_{k}\}$ an orthonormal basis of $\Hil_p$. Then $(I_{p,r} T^{-1}\psi_k)$ is an orthonormal basis of  $\Hil_r$.
		\item  If $(\psi_k)$ is a Riesz basis of $\Hil_p$  then $(I_{p,r} \psi_k)$ is a Riesz basis of  $\Hil_r$  with the same bounds. 
\end{enumerate}}
\end{cor.}

On the other hand, if $p < r$ then $\psi$ is also a sequence in $\Hil_p$. 
So, in addition to looking at $\left< f, I_{r,p} \psi_k\right>_p $ we can also look at  
$\left< f, \psi_k\right>_p $. This is done in the next section.

\subsection{Frame-related Operators for the Original Sequence $\psi$}

If we do not apply the operator we can still show
\begin{lem.}\label{lem: no chain frames}
Let $\psi = (\psi_k) \subseteq \Hil_m$ be an arbitrary sequence. Then for every $r\le p\le m$, 
\begin{equation} \label{eq:analysisrel1}
	C^r_{\psi} = C^p_{\psi}  \iota_{r,p}.
\end{equation}where $\iota_{r,p}$ is defined as in \eqref{eq: iota rp}.
\end{lem.}
\begin{proof}
We have
$$ \left< f, \psi_k\right>_r = \left< f, \iota_{p,r} \psi_k\right>_r = \left< \iota_{r,p} f, \psi_k\right>_p $$
and% for the domains we can show
\begin{equation*}
	\begin{split}
		\dom (C^p_{\psi}  \iota_{r,p})
		&= \{f\in\dom(\iota_{r,p}); \iota_{r,p}f\in\dom C^p_{\psi}\}\\
		&= \{f\in\Hil_r; (\langle \iota_{r,p}f, \psi_k\rangle_p)\in\ell^2\}\\
		&= \{f\in\Hil_r; (\langle f, \psi_k\rangle_r)\in\ell^2\}\\
		&= \dom(C^r_{\psi}).
	\end{split}
\end{equation*}
\end{proof}

In consequence, we have
\begin{lem.}\label{lem: frame propag}
If $\psi\subset\Hil_m$ and $\psi$ is a Bessel sequence for $\Hil_p$ then, for $r\leq p\leq m$, $\psi$ is a Bessel sequence for $\Hil_r$ with the same bound.
\end{lem.}
\begin{proof}
By \cite[Cor. 3.2.4 and Theor. 3.2.3]{ole1n}, 
$\psi $ is a Bessel sequence in $\Hil_p$ if and only if $\dom{D_\psi^p}=\ell^2$. By \eqref{eq:domran1} $\dom{D_\psi^p} \subseteq \dom{D_\psi^r}$, this is true if and only if $\psi$ is also a Bessel sequence for $\Hil_r$. 

Furthermore, we have that $\| C_\psi^r \| \le \| C_\psi^p \|$ by \eqref{eq:analysisrel1}. 
\end{proof} 

By \cite[Lemma 3.2]{jpaxxl09} $\psi\subset\Hil_m$ is an upper semi-frame for $\Hil_p$, $p<m$ if and only if it is a total Bessel sequence for $\Hil_p$. Then putting together Lemma \ref{lem: frame propag} and Lemma \ref{4.3.} $(i)$ we obtain 
\begin{lem.}\label{lem: up semiframe propag}Let $\psi = (\psi_{k}) \subseteq \Hil_m$ and $r\leq p\leq m$.  If $\psi$ is an upper semi-frame for $\Hil_p$, then $\psi$ is an upper semi-frame for $\Hil_r$ with the same bound.\end{lem.}

Let us now note that $\iota_{r,p}^{-1}$ is a densely defined, bijective operator, which is not bounded. It is ``never'' bounded, in the sense that if it were bounded, the involved norms would be equivalent and the whole scale of Hilbert spaces would collapse, 
however, we can use it to show 

\begin{lem.}\label{lem: frame propag 2} Let $\psi = (\psi_k) \subseteq \Hil_m$ be an arbitrary sequence. Then for every $r\le p\le m$, 
\begin{equation} \label{eq:analysisrel2}
	C^p_{\psi} = C^r_{\psi} \iota_{r,p}^{-1}\quad \text{ on }   \dom(C_\psi^p)
\end{equation} 
Therefore, ${\iota_{r,p}}_{_{|_{\dom( C^r_\psi)}}}$ is a bijective operator from $\dom(C_\psi^r)$ onto $\dom(C_\psi^p)$.  
\\\

Furthermore, if $\psi\subset\Hil_m$ and $\psi$ is complete for $\Hil_r$, then if $r\leq p\leq m$, then $\psi$ is complete for $\Hil_p$.	
\end{lem.}
\begin{proof}
The first statement is a direct consequence of \eqref{eq:analysisrel1}.

For the converse of Lemma \ref{4.3.} (i) we use \cite[Prop.4.1(g)]{xxlstoeant11}: $\psi$ is complete in $\Hil_r$ if and only if $C_\psi^r$ is injective and since $\iota_{r,p}^{-1}$ is injective, this implies that $C_\psi^p$ is injective too, which is equivalent to $\psi$ being complete in $\Hil_p$.
\end{proof} 

By Lemma \ref{lem: frame propag} we know that a Bessel sequence in $\Hil_p$ is also one in $\Hil_r$ for $r\leq p$. We can show an opposite direction for lower semi-frame. 

\begin{lem.}\label{lem: lower frame propag} Let $\psi = (\psi_k) \subseteq \Hil_m$ be a lower semi-frame for $\Hil_r$ then it is also a lower semi-frame for $\Hil_p$ for  every $r\le p\le m$ with the same lower bound. 
\end{lem.}
\begin{proof}
A sequence is a lower semi-frame if and only if the analysis operator is boundedly invertible. 
By  \eqref{eq:analysisrel2} we have that 
${C^p_{\psi}}^{-1} = \iota_{r,p} {C^r_{\psi}}^{-1}$. 

Also $\|{C^p_{\psi}}^{-1}\| \le \| {C^r_{\psi}}^{-1} \|$.
\end{proof} 

\subsubsection{Frames}
Combining those results we get
\begin{cor.}\label{lem:both frame propag} Let $\psi = (\psi_k) \subseteq \Hil_m$ be a frame for $\Hil_r$ then it is an upper semi-frame for $\Hil_p$ for  every $r\le p\le m$ and a lower semi-frame for $\Hil_q$, for every $q \le r \le m$. 
\end{cor.}

\subsubsection{Duality}

The sequence $\psi=(\psi_k)\subset\Hil_p$ is  a lower semi-frame for $\Hil_p$			if and only if \cite{CCLL} there exists a Bessel sequence $\phi=(\phi_k)\subset\Hil_p$ such that 
\begin{equation} \label{duality1}
f=\sum_k\ip{f}{\psi_k}_p\phi_k = D^p_\phi C^p_\psi (f),\quad \forall f\in\dom(C^p_\psi).
\end{equation} 

Let $\psi \subseteq \Hil_m$ be a lower semi-frame for $\Hil_r$, $r\leq m$;  then by Lemma \ref{lem: lower frame propag} it is one in $\Hil_m$. Hence there exists a dual sequence $\phi\subset\Hil_m$, which is Bessel in $\Hil_m$ and therefore also in all spaces $\Hil_r$, $r\leq m$. By assumption \eqref{duality1}  is valid on $\dom(C^m_\psi)$. 

Now let $f \in \dom(C^p_\psi)$, then by Lemma \ref{lem: lower frame propag}  \eqref{duality1} 
%$f=\sum_k\ip{f}{\psi_k}_p\phi_k$ 
converges for all $p$, with $r \le p \le m$. 
Consider $g = \io_{p,m} f$, then $g \in \dom(C^m_\psi)$.
$$ g  = \sum_k\ip{g}{\psi_k}_m \phi_k  = \sum_k\ip{ \io_{p,m} \io_{p,m}^{-1} g}{\psi_k}_m\phi_k = \sum_k\ip{\io_{p,m}^{-1} g}{\io_{m,p} \psi_k}_p \phi_k = $$
$$ =  \sum_k\ip{\io_{p,m}^{-1} g}{\psi_k}_p \phi_k = \sum_k\ip{f}{\psi_k}_p \phi_k  . $$
Furthermore
$$ f = \io_{p,m}^{-1} g = \io_{p,m}^{-1}  \sum_k\ip{f}{\psi_k}_p \phi_k. $$

So, in summary, this shows
\begin{prop.} Let $\psi = (\psi_k) \subseteq \Hil_m$ be a lower semi-frame for $\Hil_r$ then there exists a Bessel sequence $\phi \subseteq \Hil_{m}$ such that
$$ f = \io_{p,m}^{-1} \sum_k\ip{f}{\psi_k}_p\phi_k,\quad \forall f\in\dom(C^p_\psi) \subseteq \Hil_p.$$
for all $r \le p \le m$. 
\end{prop.}

\subsubsection{A negative result}

\begin{prop.}\label{lem:frame2} Let $\psi = (\psi_k) \subseteq \Hil_m$ be a frame for $\Hil_p$ and one for $\Hil_q$ , for  every $q\le p\le m$. 
Then the norms are equivalent and so $\Hil_q = \Hil_r = \Hil_p$ for  $q\le r\le p$.
\end{prop.}
\begin{proof}
By \eqref{eq:analysisrel2} $C^p_{\psi} = C^q_{\psi} \iota_{q,p}^{-1}$. As $\psi$ is a frame for $\Hil_q$ we have that $\iota_{q,p}^{-1} = ({S^{q}_\psi})^{-1}  D^q_{ \psi}  C^p_{\psi}$ and it is therefore bounded.
\end{proof}
\begin{rem.}
In particular, this means that if two Hilbert spaces, one contained into the other one, do not coincide, a sequence can {\em never} be a frame for both of them. 
\end{rem.}
\begin{rem.} For this statement it is important we have considered Hilbert spaces and (standard) frames with the sequence space $\ell^2$. If we consider Banach spaces, associated to a weighted $\ell^p$ space \cite{fe06} or Hilbert spaces with weighted $\ell^2$ sequence spaces \cite{dafora05,eh10}, this result does not hold. Quite the opposite, e.g. for localized frames \cite[Theorem 1]{xxlgro14} one can show that a frame on the pivot space $\Hil_0$ is also a frame for the other (Banach) spaces $\Hil_q$. \end{rem.}

\vspace{6pt} 

%%%%%%%%%%%%%%%%%%%%%%%%%%%%%%%%%%%%%%%%%%
%% optional
%\supplementary{The following are available online at \linksupplementary{s1}, Figure S1: title, Table S1: title, Video S1: title.}


%%%%%%%%%%%%%%%%%%%%%%%%%%%%%%%%%%%%%%%%%%

\subsection*{Acknowledgements}  G.B. acknowledges that this work has been partially supported by the 			Gruppo Nazionale per l'Analisi Matematica, la Probabilit\`{a} e le loro Applicazioni (GNAMPA) of the Istituto Nazionale di Alta Matematica (INdAM).
The work of P.B. was supported by the Innovation project of the Austrian Academy of SciencesIF\_2019\_24\_Fun  (Frames and Unbounded Operators) and the project P 34624 			(``Localized, Fusion and Tensors of Frames'' - LoFT) of the Austrian Science Fund (FWF).

\end{document}